\documentclass[conference,10pt]{IEEEtran}
\usepackage[utf8]{inputenc}
\usepackage[T1]{fontenc}
\usepackage{amsmath}
\usepackage{amsfonts}
\usepackage{array}
\usepackage{amssymb}
\usepackage{graphicx}
\usepackage{booktabs}
\usepackage{multirow}
\usepackage{geometry}
\usepackage{algorithm}
\usepackage{algpseudocode}
\author{Lingying Huang*, Xiaomeng Chen*, Wei Huo*, Jiazheng Wang, Fan Zhang, Bo Bai, Ling Shi \thanks{*These authors contributed equally to this work. L. Huang,	X. Chen, W. Huo, and L. Shi are with the Department of Electronic Engineering, Hong Kong University of Science and Technology, Clear Water Bay, Kowloon, Hong Kong, e-mails: lhuangaq@connect.ust.hk, xiaomeng.chen@connect.ust.hk, whuoaa@connect.ust.hk, eesling@ust.hk. J. Wang, F. Zhan, and B. Bai are with the Theory Lab, Huawei Hong Kong Research Centre, Hong Kong SAR, China, emails: wang.jiazheng@huawei.com, zhang.fan2@huawei.com, baibo8@huawei.com.}
\thanks{This paper has been accepted for the publication of ICARCV 2022.}}
\title{Improving Primal Heuristics for Mixed Integer Programming Problems based on Problem Reduction: A Learning-based Approach}

\IEEEoverridecommandlockouts

\begin{document}
	\maketitle
	\begin{abstract}
		In this paper, we propose a Bi-layer Predictionbased Reduction Branch (BP-RB) framework to speed up the process of finding a high-quality feasible solution for Mixed Integer Programming (MIP) problems. A graph convolutional network (GCN) is employed to predict binary variables’ values. After that, a subset of binary variables is fixed to the predicted value by a greedy method conditioned on the predicted probabilities. By exploring the logical consequences, a learning-based problem reduction method is proposed, significantly reducing the variable and constraint sizes. With the reductive sub-MIP problem, the second layer GCN framework is employed to update the prediction for the remaining binary variables’ values and to determine the selection of variables which are then used for branching to generate the Branch and Bound (B\&B) tree. Numerical examples show that our BP-RB framework speeds up the primal heuristic and finds the feasible solution with high quality.
	\end{abstract}

	\section{Introduction}
	
	Large-scale combinatorial optimization  plays a vital role in real-world problem applications,  e.g.	social network analysis \cite{rossi2014fast}, open-pit block mining \cite{jelvez2016aggregation} and scheduling medical resident training \cite{brech2019scheduling}. Owing to the large search space and NP-hardness,	such problems can be extremely challenging to solve.  Nevertheless, there exist exact solvers which are guaranteed to seek an optimal solution for combinatorial optimization problems and branch-and-bound (B\&B) is one of them. The  B\&B framework is commonly adopted to solve Mixed Integer Programming (MIP), a general formulation of  combinatorial optimization problems, and produces global optimality. In practical scenarios, it is common to repeatedly optimize homogeneous MIP instances with similar model and solution structure. Hence, it is motivated to 	adopt machine learning (ML) methods for tuning B\&B automatically  for a class of MIP instances. For example, ML has been used to improve variable \cite{khalil2016learning,nair2020solving}, node \cite{he2014learning} and cut \cite{baltean2019scoring} selection strategies. 
	
	Inside solvers, there are two prominent  sides, \textit{primal} side and \textit{dual} side. The \textit{primal} side represents the seek for good feasible solutions while the \textit{dual} side refers to the quest to prove optimality.  Though a proof of optimality is an important trait of exact solvers,  obtaining high-quality feasible solutions fast is also crucial, especially with a time limit. In real-world settings, the user will expect good feasible solutions found earlier instead of an optimal solution found many hours later. In order to ensure good primal performance, primal heuristics are equipped in modern solvers.  For instance, an open-source MIP  solver SCIP \cite{achterberg2009scip} adopts many of heuristics \cite{berthold2006primal},   including heuristics designed by experts \cite{achterberg2012rounding}, heuristics induced  by mathematical theory \cite{berthold2014rens} and meta-heuristics \cite{aarts2003local}. Extensive studies of  the computational cost of different primal heuristics and their impact on MIP solving procedure can be found in \cite{berthold2018computational,lodi2013heuristic}. 
	
	Since most primal heuristics can be computationally expensive, in the last decades, a number of works apply ML to automatically construct better primal heuristic via exploring similar structure among MIP instances, for the purpose of finding good feasible solutions quickly. Khalil et al. \cite{khalil2017learning} proposed a learning-based strategy to dynamically decided whether to run a primal heuristic at a certain node of the B\&B tree. Chmiela et al. \cite{chmiela2021learning} designed a problem-specific schedule of heuristics by learning from data describing the performance of primal heuristics. More related to our work, Shen et al. \cite{shen2021learning} trained a graph convolutional network(GCN)  to predict the optimal solution of an unseen instance and then designed a novel primal heuristic, called Probabilistic Branching technique with guided Depth-first Search (PB-DFS) which is a B\&B configuration based on a predicted solution to guide the search space of the B\&B method. 
	
	Nevertheless, the above learning-based primal heuristics are usually limited to small instances.  The huge dimensionality of a large-scale MIP instance poses a significant challenge to existing methods. In order to tackle the large-scale MIPs, we equip the PB-DFS framework in \cite{shen2021learning} with a learning-based prepossessing technique to build a Bi-layer Prediction-based Reduction Branch (BP-RB) framework, which can find high-quality feasible solutions more quickly. 
	In this paper, we focus on Integer Programming with binary variables, but it can also be adapted to other combinatorial optimization problems.  Our method works in three steps. 	First, a GCN is trained for optimal variable solution prediction based on the dataset formed by optimally solved small-scale MIPs. Then for a new MIP instance, the trained GCN model can predict for each variable its probability of taking value one or zero in the optimal solution. Second,  conditioned on the predicted probabilities, we choose  a subset of binary variables and then fix them to their predicted value. Moreover, by exploring the logical consequences of 
	fixing binary variables, we can identify the redundant constraint and fix more variables based on conflict analysis. Thus, not only the number of decision variables but  also the size of coefficients  in the constraint matrix decreases, which leads to significant reduction of the problem size.  Finally, the remaining unfixed variables define smaller `sub-MIP' and the predicted probability for the smaller sub-MIPs can be obtained by the trained GCN model. By adopting PB-DFS,  high-quality feasible solutions can be quickly found for the sub-MIP. Combining the fixed variables with the sub-MIP's feasible solution, we can obtain the feasible solution for the primal MIP instance. Our contribution are summarized as follows:	
	\begin{enumerate}
		\item We propose a bi-layer prediction-based branching framework which produces a novel B\&B configuration to seek feasible solutions quickly for large-scale combinatorial optimization problems.  
		\item We propose a learning-based problem reduction method which removes variables and constraints that are not necessary in the MIP instance. Hence, the memory required to solve the instance can be significantly reduced, speeding up the time to find quality primal solutions.  
		\item We evaluate our method on four 
		four representative NP-hard problems.  Extensive results show that our proposed method (BP-RB) can spend less time finding better primal solutions than traditional primal heuristics and PB-DFS-GCN \cite{shen2021learning}. Moreover,  we show that BP-RB can generalize very well on different problem sizes. 
	\end{enumerate}
	
	\section{Preliminaries}
	\subsection{MIP Problem}
	Consider an MIP problem following the general form:
	\begin{equation}\label{e1}
		\begin{split}
			\min \,& c^{\top}x\\
			s.t. \,& Ax\leq b,\\
			\,&	x_{j}\in\{0,1\},	\forall j\in \mathcal{B},\\
			&x_{j}\in\mathbb{Z},	\forall j\in \mathcal{Q},x_{j}\geq0,	\forall j\in \mathcal{P},		
		\end{split}
	\end{equation}
	where $ x\in\mathbb{R}^{n} $ denotes the decision variables and is partitioned into $ (\mathcal{B},\mathcal{Q},\mathcal{P}) $, with $\mathcal{B},\mathcal{Q},\mathcal{P}  $ being the index set if binary, general integer and continuous variables, respectively. We further assume that for $ j\in \mathcal{Q} $, there exist bounds $ \underline{x}_{j},\bar{x}_{j}\in\mathbb{Z} $, such that $ \underline{x}_{j}\leq x_{j}\leq \bar{x}_{j}$. Otherwise, it is easy to prove that problem \eqref{e1} is either infeasible or unbounded below. Let $ l $ represent the number of integers between $ \underline{x}_{j},\bar{x}_{j} $ including themselves and $ m=\lceil \log_{2}l\rceil $. Then solving $ x_{j} $ is equivalent to solve $ x_{j}=bx_{j,m}2^{m-1}+bx_{j,m-1}2^{m-2}+\cdots+bx_{j,1} $, where $ bx_{j,p}\in\{0,1\} $, for $ p=1,\ldots,m $. By the above analysis, any general MIP problem is equivalent to solve an MIP problem which contains binary variables and continuous variables. The main task in our paper is to predict the probability that a binary variable takes value 1 in the optimal solution.
	
	\subsection{Branch and Bound Algorithm}
	A widely-used method to produce exact solutions to MIP problems is B\&B algorithm. This method recursively builds a search tree which assigns partial integer at each node, and uses the information obtained at each node to eventually reach an optimal solution \cite{land1960automatic}. Pseudocode for the generic B\&B is given in  Algorithm \ref{agl1},
	\begin{algorithm}[htp]
		\caption{Branch and Bound} \label{agl1}
		\begin{algorithmic}[1]
			\State Set $\mathcal{L}=\mathcal{D}$ and initialize $\hat x$;
			\While{$\mathcal{L}\neq\emptyset$} 
			\State Select a subproblem $\mathcal{SP}$ from $\mathcal{L}$ to explore;
			\If{a solution $\hat x' \in\{x\in \mathcal{S}| c^{\top}x<c^{\top}\hat x\}$ can be found} 
			\State Set $\hat x =\hat x'$;
			\EndIf
			\If{$\mathcal{S}$ cannot be pruned} 
			\State Partition $\mathcal{S}$ into $\mathcal{S}_1, \mathcal{S}_2, \ldots, \mathcal{S}_r$;
			\State Insert $\mathcal{S}_1, \mathcal{S}_2, \ldots, \mathcal{S}_r$ into $\mathcal{L}$;
			\EndIf
			\State Remove $\mathcal{S}$ from $\mathcal{L}$;
			\EndWhile
			\State \Return $\hat x$.
		\end{algorithmic}
	\end{algorithm}
	where $\mathcal{D}$ is denoted as a set of valid solutions to the problem. The problem  $\mathcal{P}$ aims to find an optimal solution $x^\star \in \arg\min_{x\in\mathcal{D}} f(x)$. At each iteration, B\&B selects a new subset of the search space $\mathcal{S} \subset \mathcal{D}$ for exploration from a queue $\mathcal{L}$ of unexplored subsets. Then, if a solution $\hat x' \in \mathcal{S}$ has a better objective value than $\hat x$, the incumbent solution is updated. On the other side, the subset is pruned or fathomed if there is no solution in $\mathcal{S}$ with better objective solution than $\hat x$, i.e., $c^{\top} x\geq c^{\top}\hat x,\forall x \in \mathcal{S}$.
	
	\subsection{Primal Heuristics}
	From the above description, it is obvious that finding a good quality $ \hat{x} $ as early as possible would help to prune the search tree. The optimality of the feasible solution and the time took to find the feasible solution are two main criteria to measure the quality of primal heuristics. Some common primal heuristics are summarized by Berthold in \cite{berthold2006primal} and the comparison on the computational cost and impact on solving procedure can be found in \cite{berthold2018computational}. Here we briefly introduce some classical primal heuristics which are compared with in our simulation part.
	
	\textbf{Diving:} Diving fixes variables of fractional LP-solution to promising values and resolves LP iteratively \cite[Algorithm 1]{berthold2006primal}.
	
	\textbf{Feasibility Pump \cite{bertacco2007feasibility}:} Feasibility Pump first obtains LP-relaxation optimal solution, then rounds the integral-infeasible variable to the nearest integral. The procedure stops if the integral feasible solution is LP-feasible. Otherwise, an additional LP is solved in order to find a new LP-feasible point which is a closest to the integral feasible solution satisfying the constraint of the previous problem \cite[Algorithm 2]{berthold2006primal}. 	
	
	\textbf{Relaxation Enhanced Neighborhood Search:} This method creates a sub-MIP of the original MIP by changing the bounds of integer variables \cite[Algorithm 4]{berthold2006primal}.
	
	\textbf{Rounding:} Rounding rounds the set of fractional variables of some LP-feasible point to an integral value \cite[Algorithm 6]{berthold2006primal}.
	
	\subsection{GCN based prediction scheme}	
	\begin{figure*}[h]
		\centering
		\includegraphics[width=0.9\linewidth]{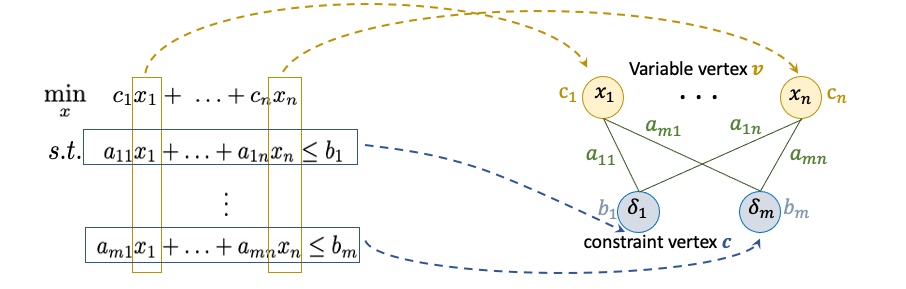}
		\caption{Bipartite graph representation of a MIP instance.}
		\label{bipartite}
	\end{figure*}

	GCN based prediction scheme  \cite{shen2021learning} predicts the probability that a binary value takes value $1$ in the optimal solution. Through treating the value prediction problem as a binary classification problem,  a dedicated graph convolutional neural network is adopted to output the probability. To model the correlations between variables from a certain type of MIP, we represent  a MIP instance as a bipartite graph, containing the objective function coefficients $\mathbf{c}$, constraint coefficient matrix $\mathbf{A}$ and right-hand-side (RHS) coefficients $\mathbf{b}$. Two sets of vertices are contained in the bipartite graph. One of them is the set of decision variable and the other is the set of constraint vertices, see Figure \ref{bipartite}. There exists an edge between a variable vertex and a constraint vertex if the corresponding variable has a nonzero coefficient in the constraint, i.e., the number of edges represents the number of non-zeros in the constraint matrix. Note that the defined bipartite graph not only captures the relation between variables and constraints, but contains the detailed coefficients in its structure as well. Particularly, the  entries in objective coefficients $\mathbf{c}$, the entries in constraint coefficients $\mathbf{b}$ and the non-zero entries of the  constraint matrix are included as scalar features of the corresponding variable vertices, constraint vertices and edges, respectively.  Besides the coefficients features, both vertices and nodes include multi-dimensional feature vectors which are extracted from the MIP, providing more information for the learning procedure. We use the code provided by Shen et al. \cite{shen2021learning}  to compute the same set of features
	using SCIP.
	
	The training dataset consists of multiple optimally-solved MIP instances from the same problem distribution and the optimal solution corresponding to each instances. One decision variable $x_i$ from a solved MIP instance is regarded as a training example and the label $y_i$ of $x_i$ is obtained from the solution value of $x_i$ in the optimal solution. Given the collected features, labels and the bipartite graphs, we can then train GCN. The cross-entropy loss function between the predicted value $\hat y_i$ of decision variables and their optimal value $y_i$ is used to train the model, which is defined as 
	\begin{equation}
		\min \quad  -\frac{1}{n}\sum\limits_{i=1}^n(y_i\times log(\hat y_i)+(1-y_i)\times log(1-\hat y_i)).
	\end{equation}

	The above GCN model ensures that the network can be applicable to MIPs of different problem sizes from the same problem distribution using the same set of parameters. That's to say, once trained on a certain type of problem dataset, the model can
	output the the probability $p_i$ that a binary variable takes value 1 in the optimal solution for an unseen instance from the same problem type. The array of predicted probabilities is referred as the \textit{probability vector}.  Based on the predicted probability, $\hat x_i$, the predicted value of variable $x_i$, equals to $1$ if  $p_i \geq 0.5$. Otherwise,  $\hat x_i=0$.

	\section{Main Results}
	\subsection{Learning-based Problem Reduction}  \label{s3.1}
	In this subsection, we  propose a novel problem reduction method based on the predicted \textit{probability vector}. The proposed problem reduction method aims to transform the large-scale problem into a smaller one so that the MIP can be solved quickly and yield high quality solution at the same time by fixing variables and removing constraints in the problem instance.
	
	\textbf{Fixing variable based on predicted probability from GCN}
	
	We sort the binary variable $x_i$ in non-increasing order of $\max (p_i, 1-p_i)$ and choose the first $\eta |\mathcal{B}|$ variables as the subset $\mathcal{F}$, where $0<\eta<1$ is a predefined value. Then, the variables with indices in $\mathcal{F}$ are fixed to their predicted values. After fixing the variables, we have an additional constant term in the objective function. Moreover, the $i-$th constraint becomes
	
	\begin{equation}
		\sum\limits_{j\in\mathcal{F}}a_{ij}\hat x_j+    \sum\limits_{s\in\mathcal{S}}a_{is} x_s \leq b_i,
	\end{equation}
	where $\mathcal{S} =\mathcal{B} \setminus \mathcal{F} $ denotes the  unfixed variables set.

	\textbf{Removing constraints}
	
	If \begin{equation}
		\sum\limits_{j\in\mathcal{F}}a_{ij}\hat x_j+    \sum\limits_{s\in\mathcal{S}}\max\{a_{is},0\}  \leq b_i,
	\end{equation} 
	then obviously the $i$-th constraint is redundant and we can remove it.

	\textbf{Fixing more variable based on logical constraints}
	
	Define $ \mathcal{S}^{+}_i=\{j|a_{ij}>0\} $ and $ \mathcal{S}^{-}_i=\{j|a_{ij}<0\} $.
	Consider a binary variable $x_k, k\in \mathcal{S}^{+}_i$, we can fix the variable $x_k$ to $0$ if
	\begin{equation}
		a_{ik}+\sum\limits_{s\in\mathcal{S}^{-}_i}a_{is}+\sum\limits_{j\in\mathcal{F}}a_{ij}\hat x_j >b_i.
	\end{equation}
	
	Next, consider a binary variable $x_k, k\in \mathcal{S}^{-}_i$, we can fix the variable $x_k$ to $1$ if
	\begin{equation}
		\sum\limits_{s\in\mathcal{S}^{-}_i\setminus \{k\}}a_{is}+\sum\limits_{j\in\mathcal{F}}a_{ij}\hat x_j >b_i.
	\end{equation}

	\subsection{Bi-layer Prediction-based Reduction Branch}
	In this subsection, we introduce our primal heuristic method. We first use GCN to predict the optimal binary solution. Then we define the score of variable $ x_{j} $ as 
	\begin{equation}
		z_{i}=\max\{p_{i},1-p_{i}\}.
	\end{equation}
	We can view this score as how certain the prediction is from the GCN model. As observed from subsection \ref{s3.1}, the score should be changed when some variables are fixed.
	
	Instead of calculating the updated score every time when fixing a variable, which would consume plenty of computational complexity, we propose a Bi-layer Prediction-based Reduction Branch (BP-RB). Pseudocode for the generic B\&B is given in  Algorithm~\ref{agl2}. 
	
	\begin{algorithm}[htb]
		\caption{Bi-layer Prediction-based Reduction Branch (BP-RB) Algorithm}\label{agl2}
		\begin{algorithmic}[1]
			\State Apply the first time GCN to predict $ z_{i} $, $ i\in\mathcal{B} $;
			\State Set threshold $ \eta $;
			\State Sort $ z_{i} $ in a descending manner, record the corresponding variable index,
			$ I\leftarrow\{i_{1},\ldots,i_{|\mathcal{B} |}\} $;
			\State Set $ \mathcal{F}\leftarrow\emptyset $;
			\While {$ j<\eta |\mathcal{B} | $}
			\State Set $ \hat{x}_{i_{j}} \leftarrow[p_{i_{j}}]$;
			\State Let $ \mathcal{F} \leftarrow\mathcal{F} \cup \{i_{j}\}$;
			\EndWhile
			\State Let $ \mathcal{S} \leftarrow \mathcal{B}\setminus\mathcal{F}$;
			\State Use learning-based problem to fix variables and remove constraints (Section \ref{s3.1}), and obtain the sub-MIP problem where its binary variables belongs to $ i\in\mathcal{B}' $;
			\State Apply the second time GCN to predict new score $ z_{i} $, $ i\in\mathcal{B}' $;
			\State Generate branch and bound tree based on $ z_{i} $ as described in Algorithm \ref{agl3};
			\State Run Depth-first Search till time limit or an optimal solution obtained, and record the optimal solution found so far as $ \hat{x} $;
			\State \Return $ \hat{x} $.			
		\end{algorithmic}
	\end{algorithm}
	
	\begin{algorithm}[htb]
		\caption{Generate B \& B Tree based on $ z $}\label{agl3}
		\begin{algorithmic}[1]
			\State Sort $ z $ in a descending manner, record the corresponding variable index,
			$ I\leftarrow\{i_{1},\ldots,i_{n}\} $;
			\State Generate the root node which represents the reduced problem, the node queue index $ \mathcal{N}\leftarrow\{0\}$;
			\For {j=1:n} \Comment{Generate branching tree according to the descending variable index}
			\State Set the right node as fixing the variable to $ x_{i_{j}}=1-[p_{i_{j}}]$;
			\State Add the right node to the node queue index $ \mathcal{N}\leftarrow\mathcal{N}\cup\{j\}$;
			\State Add constraint $ x_{i_{j}}=[p_{i_{j}}] $;			
			\EndFor
			\State Add node queue index	$ \mathcal{N}\leftarrow\mathcal{N}\cup\{n+1\}$.	
		\end{algorithmic}
	\end{algorithm}

	The main idea  is fixing the variables with high certainty. By logical constraints and other problem reduction methods described in Section \ref{s3.1}, more variables can be fixed and some constraints can be relaxed, which will lead to a sub-MIP problem with fewer integer variables and fewer constraints. As a result, the computational cost reduces significantly. In addition, the second time GCN can be viewed as taking into account the influence of the score by fixing other variables. 
	
	It is worth noting that $ |\mathcal{B}'| <| \mathcal{S} |$ due to the learning-based problem reduction method proposed in Section \ref{s3.1}. Therefore, the total number of branching variables to generate the branch tree is less than that in PB-DFS method \cite{shen2021learning}. Additionally, in the second time GCN, the number of constraints is also reduced, thus it also saves the time to solve a single node subproblem.
	
	The depth-first search method will pop out the sub-problem from the node queue to solve until time limit or an optimal solution obtained. Note that from algorithm 3, this algorithm at most will check the feasibility condition of two LP problem and solve $ n=|\mathcal{B}'| $ sub-MIP problem. Moreover, by the depth-first search method, the $ p $-th solved sub-MIP problem contains $ p $ integral variable to fix. 
	
	Different from simple diving method, this BP-RB algorithm  will explore the node most likely to contain the optimal solution first. Therefore, if the iteration time is limited, it is more likely to find a better  feasible solution (with smaller optimality gap) compared with simple diving method. By the above analysis, this method is preferable especially to find a primal heuristic solution in large-scale MIP problems.



	\section{Experimental Evaluations}
	In this section, we set up an experimental procedure followed by Shen et al.~\cite{shen2021learning} to evaluate the performance of our proposed method. First, we describe the experiment setup. Then, we illustrate the effectiveness of the proposed BP-RB by comparing it with different classical primal heuristics. Moreover, we demonstrate the efficiency of BP-RB against the full-fledged SCIP solver (SCIP-DEF), pure problem-reduction methods using the SCIP solver (ML-Split) proposed by Ding et al.~\cite{ding2020accelerating}, and PB-DFS equipped with the GCN model (PB-DFS-GCN) in Shen et al.~\cite{shen2021learning}. Further, we show the effect of different fixed proportions on the performance of BP-RB.
	
	\subsection{Setup}
	\subsubsection{Test Problems}
	We demonstrate the efficiency of our method by solving four classical combinatorial problems: Vertex Cover problem (VC), Maximum Independent Set problem (MIS), Dominant Set problem (DS), and Combinatorial Auction problem (CA). For each problem, we generate small-scale instances to train the GCN. Large-scale instances are used to evaluate the generalization of different approaches. Specifically, VC, MIS, and DS are generated based on Erd\H{o}s-R\'enyi random graphs~\cite{erdos2011evolution} with affinity $4$. The large-scale instances are formed with $3000$ variables.
	CA instances are generated by an arbitrary relationship procedure. The large-scale instances are generated with $1500$ variables and around $560$ constraints.
	
	\subsubsection{Training}
	After generating instances for each problem type, we extract $57$ statistical features of variables for each instance~\cite{shen2021learning}. Each feature is normalized and put into the GCN model for which the number of layers is set to $20$, and the dimension of the hidden vector of a binary variable is set to $32$. For a problem, we use $1000$ optimally-solved small-scale instances to train the GCN.
	
	\subsubsection{Evaluation of Solution Methods}
	We compare the performance of our BP-RB algorithm with traditional primal heuristic methods and PB-DFS-GCN by Shen et al.~\cite{shen2021learning} on large-scale problem instances. 
	
	First, we compare BP-RB with traditional primal heuristics without requiring a feasible solution on large-scale problem instances. Specifically, four heuristics are chosen as baselines: Relaxation Enhanced Neighborhood Search (RENS), Feasibility Pump (FSP), Diving, and Roundings.
	
	Second, BP-RB is compared with a full-fledged SCIP solver in which all heuristics are enabled, SCIP solver with only problem reduction, and SCIP solver with only PB-DFS-GCN.
	
	For all the above comparisons, $40\%$ variables are fixed in BP-RB before features of sub-MIP problems are put into the GCN. Further, we compare the performance of BP-RB under different portions of fixed variables.

	\subsubsection{Experimental Environment}
	All experiments are conducted on a cluster of $64$ Intel $3.40$ GHz CPUs and 16GB RAM. Our method is implemented via C-api provided by SCIP, version $6.0.1$. The GCN is implemented by the Tensorflow package.
	\subsection{Evaluation Results for Finding Primal Solutions}
	Comparisons of our method with effective heuristic methods in SCIP solver are illustrated in Table~\ref{tab:primal_heuristic}. Note that the PB-DFS-GCN approach is applied only at the root node. It stops running once finding the first feasible. Our method first reduces the number of variables and constraints in an original problem to yield a subproblem, and then capitalize on PD-DFS at the root node to seek the first feasible solution. For each type of MIP problem, we run $30$ different instances for a heuristic, and each run terminates with a cutoff time of $50$ seconds. \textbf{\# Instances no feasible solution} shows the number of instances that the corresponding heuristic does not find any feasible solution. Heuristics that does not find any primal feasible solution for all $30$ instances is omitted in the table. Results in other columns are geometric mean shifted by one averaged over solved instances. Here, we not only consider the best solutions but also look at the best heuristic solution objective found by a heuristic. For best heuristic solution time, best solution time, and heuristic total time, their time is the summation of problem reduction time and PB-DFS-GCN running time.
	
	As shown in Table~\ref{tab:primal_heuristic}, overall, BP-RB can spend less time finding better primal solutions than traditional primal heuristics and PB-DFS-GCN not only on VC, DS, MIS but also on CA which is a kind of problem not formulated on graphs. 
	Note that PB-DFS-GCN~\cite{shen2021learning} is less competitive on CA in terms of both objective solution and running time, but BP-RB performs better due to problem reduction and probability adjustment in subproblems. These results show that BP-RB can generalize very well on problem size and problems not related to graphs. 
	\begin{table*}[htbp] 
		\centering  
		\caption{\label{tab:primal_heuristic}Comparison of BP-RB with primal heuristics and PB-DFS-GCN}   
		\begin{tabular}{m{0.08\textwidth}<{\centering} m{0.15\textwidth}<{\centering} m{0.1\textwidth}<{\centering}m{0.1\textwidth}<{\centering} m{0.1\textwidth}<{\centering} m{0.08\textwidth}<{\centering} m{0.1\textwidth}<{\centering}m{0.08\textwidth}<{\centering}}    
			\toprule 
			Problem & Heuristic &Best Heuristic Solution Objective & Best Heuristic Solution Time & Best Solution Objective & Best Solution Time & \# Instances no feasible solution & Heuristic Total Time\\    
			\midrule
			\multirow{4}{*}{VC (Min.)} 
			&Roundings & 1813.3 & 32.0 & 1774.0& 38.6& 0 & 6.6\\ 
			&Feaspump& 2137.3& 2.5& 1767.6& 43.0& 17&1.7 \\
			&PB-DFS-GCN & 1628.2& 8.4 &1628.2& 8.4& 0 & 8.4\\   
			&BP-RB & \bf{1627.8} & \bf{0.7} & \bf{1626.8} & \bf{0.9} &0&\bf{0.7} \\ 
			\midrule
			\multirow{4}{*}{DS (Min.)} 
			&Roundings & 622.5 & 21.1 & 622.5 & 21.2 & 0 & \bf{0.5}\\ 
			&Feaspump& 325.1 & 8.3 & 325.1 & \bf{8.3} & 0 & 0.8 \\
			&PB-DFS-GCN & \bf{318.4} & 18.9 & \bf{318.4} & 18.9 & 0 & 18.9 \\   
			&BP-RB & 319.8 & \bf{5.2} & 319.3 & 9.8& 0& 5.2 \\ 
			\midrule
			\multirow{5}{*}{MIS (Max.)}
			&Roundings & 995.6 & 36.0 & 1113.6 & 37.2 & 0 & 7.2\\ 
			&Feaspump& 845.1  & 2.8 & 1022.6 & 12.3 & 1 & 1.9 \\
			&Diving & 829.6 & 46.0 & 1228.1 & 42.7 & 22 & \bf{0.1}\\  
			&PB-DFS-GCN & \bf{1369.0} & 8.4& \bf{1369.0} & 8.4 & 0 & 8.5\\   
			&BP-RB & 1367.4 & \bf{2.2} & 1367.7 & \bf{2.2} & 0& 2.2 \\ 
			\midrule
			\multirow{5}{*}{CA (Max.)}&RENS& 3425.9 & 6.3 & 3674.3 & 35.5 & 0 & 3.9 \\   
			&Roundings & 3206.0 & 36.2 & 3688.8 & 32.4 & 0 & \bf{0.3}\\ 
			&Diving & \bf{3481.6} & 10.3 & 3698.2 & 36.6 & 0 & 0.7\\  
			&PB-DFS-GCN & 3341.5& 10.5& 3664.9 & 37.6 & 0 & 10.6\\   
			&BP-RB & 3368.7 & \bf{2.0} & \bf{3718.3} & \bf{13.6}& 0& 2.0 \\
			\bottomrule 
		\end{tabular}  
	\end{table*}
	
	We further compare BP-RB with the use of full-fledge SCIP, SCIP equipped with only problem reduction, and SCIP with only PB-DFS-GCN. Here, we do not limit the running time, but let each heuristic help the SCIP to solve optimization problems. The detailed solving statistics results are given in Table~\ref{tab:SCIP_fledge}. We observe that our method significantly outperforms other methods on VC, MIS, and CA. Finding a good primal feasible solution early helps BP-RB solve large-scale MIP problems faster.
	\begin{table}[htbp] 
		\centering  
		\caption{\label{tab:SCIP_fledge}BP-RB compared with SCIP-DEF, ML-Split and BP-RB}   
		\begin{tabular}{m{0.15\linewidth}<{\centering} m{0.22\linewidth}<{\centering} m{0.12\linewidth}<{\centering} m{0.12\linewidth}<{\centering} m{0.12\linewidth}<{\centering}}     
			\toprule 
			Problem & Method & Best Solution Objective &Best Solution Time & Optimally Gap (\%) \\  
			\midrule
			\multirow{4}{*}{VC (Min.)}& SCIP-DEF & 1634.4 & 508.3 & 3.8 \\ 
			& ML-Split & 1630.8 & 120.2 & 3.6\\  
			& PB-DFS-GCN & 1628.2 & 3.7 & 3.3 \\
			& BP-RB & \bf{1626.8} & \bf{0.9} & \bf{0.0} \\
			\midrule
			\multirow{4}{*}{DS (Min.)}& SCIP-DEF & 315.3 & 388.9 & 2.9 \\ 
			& ML-Split & 316.3 & 274.8 & 3.2 \\  
			&PB-DFS-GCN & 315.7 & \bf{243.4} & 2.9\\   
			&BP-RB & \bf{315.2} & 537.6 & 2.3\\  
			\midrule
			\multirow{4}{*}{MIS (Max.)}& SCIP-DEF &  1362.4 & 591.4& 4.8  \\ 
			& ML-Split & 1359.3 & 121.9  & 5.1 \\  
			& PB-DFS-GCN & \bf{1369.0} & 3.7& 4.2 \\
			& BP-RB & 1367.7 & \bf{2.2} & \bf{0.0} \\ 
			\midrule
			\multirow{4}{*}{CA (Max.)}& SCIP-DEF & 3825.6 & 794.4 & 11.3 \\ 
			& ML-Split & 3825.8 & 685.8 & 11.7  \\  
			& PB-DFS-GCN & 3825.6 & 794.4 & 11.3 \\
			& BP-RB & \bf{3896.8} & \bf{621.8} & \bf{7.3} \\
			
			\bottomrule   
		\end{tabular}  
	\end{table}	
	
	\begin{table}[htbp]   
		\centering
		\caption{\label{tab:portion}Comparison of different number of fixed variables}   
		\begin{tabular}{m{0.15\linewidth}<{\centering} m{0.15\linewidth}<{\centering} m{0.15\linewidth}<{\centering} m{0.15\linewidth}<{\centering} m{0.15\linewidth}<{\centering}}     
			\toprule 
			Problem & Fixed Portion & Best Heuristic Solution Objective & Best Heuristic Solution Time & Best Solution Objective \\    
			\midrule
			\multirow{4}{*}{VC (Min.)}& $20\%$ & 1628.5 & 0.9 & 1627.3 \\ 
			& $40\%$ & \bf{1627.8} & \bf{0.7} & \bf{1626.8}\\  
			& $60\%$ & 1633.0 & 0.9 & 1633.0 \\
			& $80\%$ & 1634.1 & 1.0 & 1634.1 \\
			\midrule
			\multirow{4}{*}{DS (Min.)}& $20\%$ & \bf{318.8} & 5.6 & \bf{318.5} \\ 
			& $40\%$ & 319.8 & 5.2 & 319.3 \\  
			& $60\%$ & 322.3 & 3.0 & 320.9 \\
			& $80\%$ & 330.6 & \bf{0.7} & 323.4 \\  
			\midrule
			\multirow{4}{*}{MIS (Max.)}& $20\%$ & \bf{1367.7} & \bf{2.14} & \bf{1370.1} \\ 
			& $40\%$ & 1367.4 & 2.2 & 1367.7 \\  
			& $60\%$ & 1366.4 & 3.23 & 1366.6 \\
			& $80\%$ & 1340.9 &  3.9 & 1340.9 \\ 
			\midrule
			\multirow{4}{*}{CA (Max.)}& $20\%$ & 3325.9 & 2.4 & 3708.3 \\ 
			& $40\%$ & 3368.7 & 2.0 & 3718.3 \\  
			& $60\%$ & \bf{3416.5} & 1.1 & \bf{3854.7} \\
			& $80\%$ & 3411.06 & \bf{0.3} & 3780.66 \\
			
			\bottomrule   
		\end{tabular}  
	\end{table}
	
	Moreover, comparisons of the performance of the proposed BP-RB under different portions of fixed variables are given in Table~\ref{tab:portion}. We can observe that the best threshold is different for different types of problems. For example, VC performs best when we fix $40\%$ variables to reduce the problem, while MIS performs best when $60\%$ variables are fixed. Further, there is a trade-off between the quality of primal feasible solution and the time to find the primal feasible solution. For instance, for DS problems, increasing the number of fixed variables significantly reduces the best heuristic solution time, but leads the best heuristic solution far away from the optimal solution.

	\section{Conclusions}
	In this paper, we propose a Bi-layer Prediction-based Reduction Branch (BP-RB) framework in order to accelerate the process of finding a high-quality feasible solution. A graph constitutional network (GCN) is employed to predict the binary variables' values. The GCN framework has generalization of different variable dimensions, which makes it possible to train the network with small-scale problems. In addition, conditioned on the predicted probabilities, a subset of binary variables are fixed to the predicted value by greedy method. By exploring the logical consequences, a learning-based problem reduction method is proposed. This method significantly reduces the variable and constraint size. With the produced reductive sub-MIP problem, the second layer GCN framework is employed to update the prediction for the leaving binary variables' values  which serves as the variable selection for branching criterion. Numerical examples show that our BP-RB framework speeds up the primal heuristic with high quality. We also show the comparison on the objective quality and time between different fixed portion. However, the optimal fixed portion changes for different problems. We will look for an adaptive fixing portion for different problems considering the time and quality constraint in future work.	
	
	
	\bibliographystyle{ieeetr}
	\bibliography{ref}

\end{document}